\documentclass{amsart}
\usepackage{ upgreek }
\usepackage{latexsym,dsfont}
\usepackage{amsmath,amssymb,amsfonts}
\usepackage[utf8]{inputenc}
\usepackage[colorlinks=true, urlcolor=blue, citecolor=red]{hyperref}
\usepackage{cleveref}
\usepackage{latexsym,dsfont}
\usepackage{algorithm}
\usepackage{algpseudocode}

\usepackage{subcaption}[skip=1pt]

\newcommand{\cH}{\mathcal{H}}
\usepackage{enumitem}

\DeclareMathOperator*{\argmax}{arg\,max}

\newcommand{\bfa}{\boldsymbol{a}}
\newcommand{\bfb}{\boldsymbol{b}}
\newcommand{\bfc}{\boldsymbol{c}}

\newcommand{\bfh}{\boldsymbol{h}}

\newcommand{\bfp}{\boldsymbol{p}}

\newcommand{\bfs}{\boldsymbol{s}}

\newcommand{\bfu}{\boldsymbol{u}}

\newcommand{\bfw}{\boldsymbol{w}}
\newcommand{\bfx}{\boldsymbol{x}}
\newcommand{\bfy}{\boldsymbol{y}}
\newcommand{\bfz}{\boldsymbol{z}}

\newcommand{\bfK}{\boldsymbol{K}}

\newcommand{\bfQ}{\boldsymbol{Q}}

\newcommand{\dd}{\rm d}
\newcommand{\Omt}{\Omega_{\rm targ}}

\newcommand{\pith}{\pi_{\bfth}}

\newcommand{\bfth}{\boldsymbol{\theta}}

\newcommand{\bfell}{\boldsymbol{\ell}}

\newcommand{\R}{\ensuremath{\mathds{R}}}
\newcommand{\E}{\ensuremath{\mathds{E}}}

\newtheorem{theorem}{Theorem}[section]

\theoremstyle{definition}

\theoremstyle{remark}
\newtheorem{remark}[theorem]{Remark}
\newtheorem{assumption}[theorem]{Assumption}

\usepackage[textsize=small]{todonotes}
\usepackage{lineno}
\usepackage[dvips,margin=1.5in]{geometry}
\numberwithin{equation}{section}

\newcommand{\fig}[1]{\Cref{#1}}

\title[Neural Network Approaches for Parameterized Optimal Control]{ Neural Network Approaches for Parameterized Optimal Control}

\author[D. Verma]{Deepanshu Verma$^\ast$}

\author[N. Winovich]{Nick Winovich$^\dag$}

\author[L. Ruthotto]{Lars Ruthotto$^\ddag$}

\author[B. van Bloemen Waanders]{Bart van Bloemen Waanders$^\dag$}
\address{$^\ast$Department of Mathematics, Emory University, Atlanta, GA}
\email{dverma4@emory.edu}
\address{$^\dag$Departments of Mathematics and Computer Science, Emory University, Atlanta, GA}
\email{lruthotto@emory.edu}
\address{$^{\ddag}$Sandia National Laboratories, Albuquerque, NM}
\email{\{nwinovi, bartv\}@sandia.gov}

\date{}

\begin{document}

\begin{abstract}
  We consider numerical approaches for deterministic, finite-dimensional optimal control problems whose dynamics depend on unknown or uncertain parameters.
    We seek to amortize the solution over a set of relevant parameters in an offline stage to enable rapid decision-making and be able to react to changes in the parameter in the online stage.
    To tackle the curse of dimensionality arising when the state and/or parameter are high-dimensional, we represent the policy using neural networks.
    We compare two training paradigms:
    First, our model-based approach leverages the dynamics and definition of the objective function to learn the value function of the parameterized optimal control problem and obtain the policy using a feedback form.
    Second, we use actor-critic reinforcement learning to approximate the policy in a data-driven way.
    Using an example involving a two-dimensional convection-diffusion equation, which features high-dimensional state and parameter spaces, we investigate the accuracy and efficiency of both training paradigms.
    While both paradigms lead to a reasonable approximation of the policy, the model-based approach is more accurate and considerably reduces the number of PDE solves.
    
\end{abstract}
\maketitle
\section{Introduction}
We are interested in deterministic, finite-time optimal control problems in which the dynamical system depends on parameters that are unknown or uncertain.
Such problems arise in decision-making for complex systems and are prevalent across application domains including physics, engineering, finance, robotics, economics, environmental sciences, and traffic flow \cite{traffic-flow-2010,Trltzsch2010OptimalCO,rozza2012reduction,Ballarin2017NumericalMO,strazzullo2018model,xu2021machine}.

We aim to amortize the solution across multiple parameter uncertainties and obtain a policy that maps the current state and parameter to approximate optimal actions.  While this introduces significant offline costs, the policy is fast to evaluate online so that rapid decision-making is enabled. Moreover, our capability supports changes in the parameters over time.

Since an effective policy must depend on the current state of the system and the value of the parameter, amortizing control policies are affected by the curse of dimensionality.
When the sum of the dimensions of the state or parameter exceeds approximately four, traditional solution approaches become ineffective.
Tackling the curse of dimensionality leads to two critical challenges:
First, one needs to determine policies using a high-dimensional function approximator.
Motivated by their universal approximation properties, we consider neural networks.
Second, identifying the relevant parts of the state space is critical to avoid the infeasibility of globally sampling high-dimensional spaces. 
Therefore, it is common to limit the search to parts of the state space that are likely to be visited when following optimal policies~\cite{Li-Verma-Ruthotto-2022, sutton2018reinforcement}.

Given a class of function approximators, we focus primarily on numerical methods for learning effective policies and compare a model-based and a data-driven approach. 
Our model-based approach approximates the value function of the control problem from which the policy is obtained using a feedback form.
The feedback form requires a model of the dynamical system and knowledge of the objective functional.
This enables the incorporation of the physics of the system and include known properties of the value function, via the Hamilton-Jacobi-Bellman (HJB) equations, into the training problem.
Our data-driven approach uses the framework of actor-critic reinforcement learning to decouple the networks for the policy (i.e., the actor) and the value function (i.e., the critic) and trains both networks purely from observations.
Both approaches are intimately related through the underlying control theory, especially Pontryagin's maximum principle theory and HJB equations.

To compare the accuracy and efficiency of both approaches, we implement a model problem motivated by control of contaminant transport.  Here, the state is governed by the two-dimensional advection-diffusion partial differential equation (PDE).  The objective is to prevent a contaminant originating at a source location from reaching a target region by controlling the location of a sink.  Uncertain parameters in this case consist of different scenarios for the source location and the velocity field.  Since both the state and the parameter are infinite-dimensional, the problem is well-suited to create high-dimensional problem instances through discretization. Both the state and the parameters are infinite-dimensional and therefore motivate high-dimensional problems after discretization.

Our main contributions and observations are summarized as follows:
\begin{itemize}
    \item We extend the optimal control approaches from~\cite{onken2020neural} to parametric control problems, which leads to a considerable increase in the dimensionality of the problem. 
    \item We provide a direct comparison of the OC and RL approaches and including resective sub-optimalities. For the parameterized OC approach, this experiment increases the dimension from around 100 in~\cite{onken2020neural} to about 1000 (see \cref{sec:num-exp}). 
    \item We observe in our experiments that the OC approach achieves more accurate policies with fewer simulations (20x less) of the dynamical system. 
\end{itemize}

The remainder of the paper is organized as follows: \Cref{sec:OCF} introduces the parametric OC problem and provides a brief background on OC theory that informs our approaches. \Cref{sec:related} discusses related approaches for learning amortized policies. \Cref{sec:HJB-approach} presents our neural network approach for computing the value function. \Cref{sec:RL} revisits key RL concepts, laying the groundwork for employing Proximal Policy Optimization (PPO) and Twin Delayed Deep Deterministic Policy Gradients (TD3) algorithms in policy approximation. \Cref{sec:num-exp} conducts a comprehensive numerical assessment of both the neural-HJB and RL approaches, employing PDE-based dynamics. Finally, \Cref{subsec:results} provides a thorough discussion of results, highlighting the strengths and limitations of the proposed methodologies.

\section{Parameterized Optimal Control Problem}\label{sec:OCF}
 We consider a family of deterministic, finite horizon optimal control problems of the form 
 \begin{subequations}\label{eq:OCP}
\begin{align}\label{eq:ctrl-obj}
\min_{\bfu \in U} J(t,\bfx,\bfu;\bfy) := \int_{t}^{T} L(s,\bfz(s),\bfu(s)) \dd s+  \,G(\bfz(T)),
\end{align}
subject to a nonlinear dynamical constraint
\begin{equation}\label{eq:sys_dyn}
\begin{aligned}
\dd_s \bfz(s) = f (s,\bfz(s);\bfy) + g(s,\bfz(s);\bfy)\; \bfu(s), \quad s\in (t,T];\quad
\bfz(t) = \bfx,
\end{aligned}
\end{equation}
\end{subequations}
where, given a parameter $\bfy \in \eta$, the functions $f:[0,T]\times \R^d \to \R^d$ and $g:[0,T]\times \R^d \to \R^{d\times a}$ model the evolution of the states $\bfz:[0,T]\to \R^d$ in response to the controls $\bfu:[0,T]\to U \subset \R^a$. The initial state of the system is denoted as $\bfx \in \R^d$, and the initial time is represented by $t \in [0,T]$. The function $L:[0,T]\times \R^d \times \R^a\to \R$ represents the
running cost, and $G:\R^d \to \R$ represents the terminal cost in the control objective \cref{eq:ctrl-obj}.

We assume that $f, L, G, U$ are sufficiently regular (see \cite[Chapter~2]{yong_zhou-1999} for a list of assumptions)

We aim to determine an optimal control (policy) amortized over the parameters $\bfy \in \eta$ that minimizes the overall cost. In particular, given a specific parameter $\bfy$, we seek to find the value function
\begin{equation}\label{eq:val_fn}
    \Phi(t,\bfx;\bfy) = \inf_{\bfu \in U} J(t,\bfz,\bfu;\bfy) \quad \text{s.t. \quad \cref{eq:sys_dyn}} ,
\end{equation}
which represents the minimum cost-to-go for any initial state $(t,\bfx)$. The control $\bfu^*$ for this minimum value is referred to as optimal control, and the corresponding trajectory $\bfz^*$ is called an optimal trajectory. 

The problem \cref{eq:OCP} is typically tackled using local solution methods. These methodologies involve minimizing the cost functional by eliminating the constraints, all while considering a fixed parameter $\bfy$. However, these localized strategies necessitate reevaluating the solution each time the $\bfy$ changes, rendering the previously computed solution obsolete. This clearly highlights the need for devising a control policy that is amortized over the parameters $\bfy \in \eta$. 

Subsequently, we will see that the value function $\Phi$ contains complete information about the optimal control, and we can straightforwardly deduce $\bfu^*$ from $\Phi$. Notably, this connection is not explicitly utilized by RL algorithms during policy approximation. This distinction between the HJB approach and RL is significant, as leveraging the underlying physics of the problem can lead to improved results as presented in \cref{sec:num-exp}.

We now introduce the Hamiltonian, $H :[0,T] \times \R^d \times \R^d \to \R \cup \{\infty\}$ of the system \cref{eq:OCP}, which plays a pivotal role in both OC theory and our numerical approach:
\begin{equation}\label{eq:H}
		H(s,\bfz,\bfp;\bfy)=\sup_{\bfu \in U} \cH(s,\bfz,\bfp,\bfu;\bfy)
\end{equation}
Here $\bfp:[t_0,T]\to \R^d$ is referred to as the adjoint or costates of the system, and
\[ \cH(s,\bfz,\bfp,\bfu;\bfy) =  \bfp \cdot \left( f(s,\bfz;\bfy) + g(s,\bfz;\bfy)\; \bfu \right)-L(s,\bfz,\bfu).  \]
The Pontryagin maximum principle provides a set of first-order necessary conditions and states that, see \cite[Theorem~I.6.3]{flemingsoner06}, the adjoint, $\bfp$, satisfies the adjoint equation. Furthermore, when the value function $\Phi$ is differentiable, the adjoint $\bfp$ and the optimal control $\bfu^*$ can be derived from $\Phi$ using \cite[Theorem I.6.2]{flemingsoner06}: 
for $s\in[t,T]$
\begin{equation}\label{eq:gradPhi}
\bfp(s)=\nabla_{\bfz} \Phi \big(s,\bfz^*(s); \bfy \big), \quad \text{and}
\end{equation}
\begin{equation}\label{eq:opt_control_nec}
\begin{split}
\bfu^*(s, \bfz^*(s),\nabla_{\bfz} \Phi \big(s,\bfz^*(s));\bfy)\in \argmax_{\bfu\in U} \mathcal{H}\big( s,\bfz^*(s),\nabla_{\bfz} \Phi \big(s,\bfz^*(s); \bfy\big),\bfu(s);\bfy).
\end{split}
\end{equation}
The \cref{eq:opt_control_nec} characterizes optimal control in feedback form, and states that it can be readily computed at any space-time pair, for a given $\bfy$, when the value function $\Phi$ and $\nabla_{\bfz} \Phi$ are available. For a thorough exploration of the aforementioned relationships, we refer the readers to \cite{flemingsoner06, Li-Verma-Ruthotto-2022, onken2020neural}.

The model-based approach exploits this fundamental relationship to approximate a $\Phi$ amortized across diverse parameters $\bfy$. This ensures a readily accessible approximate optimal control for each parameter setting. The data-driven approach, on the other hand, is agnostic to this relationship and hence suffers in terms of accuracy.

\section{Related Work} \label{sec:related}
This section reviews some closely related approaches for control/policy training. 
Existing approaches can roughly be grouped by the type of function approximators used and whether they primarily use the model or data during training. Specifically focusing on NNs as function approximators, we distinguish between model-based and data-driven training algorithms.
 
Model-based approaches, like the one presented by Kang \textit{et al.} \cite{kang2019algorithms} adopt a supervised learning approach to learn a closed-form value function. The physics-informed neural networks (PINNs) method \cite{lu2021physics,mowlavi2022optimal} integrates the cost functional into the standard PINNs loss for solving OC problems. Other techniques, such as deep-learning-based surrogate models \cite{xu2021machine,lye2021iterative} and operator learning methods \cite{wang2021fast,hwang2022solving}, focus on achieving fast OC solution inference without intensive computations. In \cite{demo2021extended}, an extended PINN approach utilizes Karush–Kuhn–Tucker (KKT) conditions as the loss function. These approaches primarily tackle OC problems with partial differential equations (PDEs) as constraints.

A recent development is the adjoint-oriented neural network method proposed in \cite{Yin2023AONNAA} for parametric OC problems. This method learns the states, adjoints, and control of the system through three different neural networks using a direct-adjoint looping method-type loss. While capable of handling parametric problems simultaneously with random sampling instead of spatial domain discretization, it faces challenges in the presence of changes in initial states or disturbances when following the optimal solution. Moreover, training separate networks introduces computational complexity, limiting the applicability of this approach to less complex problems.

Reinforcement learning (RL) has emerged as a potent data-driven approach and shown success in diverse applications such as robotics, autonomous driving, games, and recommender systems \cite{Mnih2013PlayingAW,Silver2016MasteringTG,Wu2017LearningTS,Silver2018AGR,Chen2018TopKOC,Hwangbo2019LearningAA,Berner2019Dota2W}. Its non-intrusive and gradient-free nature make RL suitable for problems lacking well-defined models. However, RL often requires many samples, which are computationally expensive in our setting. The process of hyperparameter tuning in RL algorithms is laborious and computationally expensive, hindering their applicability for parametric dynamics \cite{nikishin2018improving,henderson2018deep,rlblogpost,deepRL-chapter7-2020,amin2021survey}. To overcome these limitations, there is a growing interest in effectively integrating system knowledge to overcome the constraints of data-driven methods for handling complex dynamics. In \cite{Zhou_2021,zhou2023policy}, the authors investigate solving OC problems using actor-critic and policy gradient type methods within the realm of RL. While successful for general OC problems, these methods have limited direct applicability to parametric problems. The main challenge arises from the added complexity introduced by parameters, in addition to the already expensive task of evolving the dynamics over time.

Works closely related to our OC approach include \cite{onken2020neural,kunisch2023optimal,Li-Verma-Ruthotto-2022}, which employ NNs to parameterize the value function and penalize the HJB equation satisfied by the value function.

\section{Model-Based Approach}\label{sec:HJB-approach} 
In this section, we present our model driven approach, inspired by \cite{onken2020neural}, for approximating the amortized value function of the parameterized OC problem \cref{eq:OCP}. The central concept revolves around utilizing an efficient function approximator to model the amortized value function $\Phi$ in \cref{eq:val_fn} and subsequently computing the control using the feedback form given by \cref{eq:opt_control_nec}. The learning adopts an unsupervised approach, akin to RL. The loss function consists of the sum of the expected cost over different parameters, $\bfy\in \eta$, driving the trajectories and penalty terms that enforce the HJB equations along the trajectories and at the final-time.

\subsection{Learning Problem} \label{sec:HJB-train}
Let $\Phi_{\bfth}$ denote the approximation of the value function $\Phi$ with parameters $\bfth$. Ideally, we aim for a choice of $\bfth$ where $\Phi_{\bfth}$ matches the value function of the corresponding control problem for every given $\bfy$. 
However, this problem suffers from the curse of dimensionality so we resort to an approach which enforces this property in a subset of space-time domain.

To learn the parameters in an unsupervised way (that is, no apriori data for $\Phi$ and optimal control trajectories), we approximately solve the minimization problem
\begin{equation}\label{eq:train_loss}
\begin{aligned}
    	&\min_{\bfth} \E_{\bfy \sim \eta}   \left\{
       J(t,\bfz,\bfu;\bfy) +  P_{\rm HJB, \, \bfy}(\bfz) \right\} \\
\text{s.t.}  \quad &\dd_s \bfz = f (s,\bfz(s);\bfy) + g(s,\bfz(s);\bfy)\; \bfu(s), \quad s\in (t,T];\quad
\bfz(t) = \bfx.
      \end{aligned} 
\end{equation}
Here $P_{\rm HJB, \, \bfy}$ penalizes deviations from the HJB PDE satisfied by the value function for each parameter $\bfy$. In the learning process, the optimal control is computed through the feedback form provided in \cref{eq:opt_control_nec},
hence, incorporating the model information.

Exploiting the fact that the value function satisfies the HJB PDE (see \cite[Theorems I.5.1, I.6.1]{flemingsoner06}), we guide the approximation of $\Phi$ using the penalty term $P_{\rm HJB, \, \bfy}$, defined as:
\begin{equation}\label{eq:p_hjb}
\begin{split}
     P_{\rm HJB, \, \bfy}(\bfz) = &\beta_1 \int_t^T |H(s,\bfz,\nabla \Phi_{\bfth}(s,\bfz(s);\bfy);\bfy) - \partial_s \Phi_{\bfth}(s,\bfz(s);\bfy)| \dd s \\ 
     &+\beta_2  |G(\bfz(T)) - \Phi_{\bfth}(T, \bfz(T);\bfy)| +  \beta_3 |\nabla_{\bfz} G(\bfz(T)) - \nabla_{\bfz} \Phi_{\bfth}(T,\bfz(T);\bfy)|,
\end{split}
\end{equation}
where the relative influence of each term is controlled by the components of $\beta=(\beta_1,\beta_2,\beta_3) \in \R^3_{+}$.

We further make the following assumption for our model driven approach.
\begin{assumption}
    There exists a closed form solution to \cref{eq:opt_control_nec}, which allow us to write the control $\bfu$ as a function of $\nabla_{\bfz} \Phi$, explicitly, i.e.
    \begin{equation}\label{eq:opt_cont_nec_close}
        \bfu^*(s, \bfz^*(s),\nabla_{\bfz} \Phi \big(s,\bfz^*(s);\bfy);\bfy) = \argmax_{\bfu\in U} \mathcal{H}\big( s,\bfz^*(s),\nabla_{\bfz} \Phi \big(s,\bfz^*(s); \bfy\big),\bfu(s);\bfy),
    \end{equation}
\end{assumption}
A closed-form solution for the optimal control exists in a wide variety of OC problems~\cite{lopez2019solutions,bansal2021deepreach,kunisch2020semiglobal}.  While we do not explicitly demonstrate it in this work, it is worth noting that this assumption can be relaxed to include implicitly defined functions as long as they can be efficiently obtained. This flexibility allows for the modeling of more general convex running costs and enhances the applicability of our approach to a wider range of problems.

The \cref{eq:opt_cont_nec_close,eq:train_loss} provide a framework for obtaining the optimal control and trajectory based on the value function $\Phi$, subject to certain smoothness assumptions, by outlining the necessary conditions for optimality. Once a good approximation of the value function $\Phi$ is computed, this framework can be applied to any initial data and parameters, also allowing for adaptability to perturbations in the system. 
\subsection{Function Value Approximation} \label{sec:HJB-NN-arch}
In principle, one could employ any high-dimensional function approximator to parameterize the value function. Yet, due to the universal approximation properties inherent in NNs, we opt for using an NN to parameterize the value function.
Designing an effective neural network architecture is essential for various learning tasks, and it remains an active area of research. In our approach, we treat this as a modular component, providing flexibility. Our framework can seamlessly integrate with any scalar-valued neural network that accepts inputs in $\R^{d+1}$ and possesses at least one continuous derivatives concerning its first $d+1$ inputs, to allow computations of $\nabla \Phi$.

In our experiments, we use a residual neural network given by
\begin{equation}
\label{eq:NNArchitecture}
	\Phi_{\bfth}(\bfh_0) = \bfw^\top {\rm NN}(\bfh_0;\bfth_{\rm NN}),
\end{equation}
with trainable weights $\bfth$ containing $\bfw\,\,{\in}\,\,\R^m$ and $\bfth_{\rm NN}\,\,{\in}\,\,\R^p$.
Here the inputs $\bfh_0=(t,\bfz(t);\bfy)\in \R^{d+1}$ correspond to time-space, and  ${\rm NN}(\bfh_0;\bfth_{\rm NN}) \colon \R^{d+1} \to \R^m$ is a residual neural network (ResNet)~\cite{He_2016_CVPR}
\begin{equation} \label{eq:ResNet}
    \begin{split}
    \bfh_1 & = \rm \sigma(\bfK_0 \bfh_0 + \bfb_0) \\ 
    \bfh_{i+2} & = \bfh_{i+1} + \rm \sigma(\bfK_{i+1} \bfh_{i+1} + \bfb_{i+1}),\quad 0\le i \le M-2 \\
    {\rm NN}(\bfh_0;\bfth_{\rm NN}) & = \bfh_{M} + \rm \sigma(\bfK_M \bfh_{M} + \bfb_M), 
    \end{split}
\end{equation}	
with neural network weights $\bfth_{\rm NN}{=}(\bfK_0,\hdots,\bfK_M, \bfb_0,\hdots,\bfb_M)$ where $\bfb_i\in \R^{m}\; \forall i$, $\bfK_0 \in \R^{m \times (d+1)}$, and $\{ \bfK_1,\hdots,\bfK_M \}\in \R^{m \times m}$ with $M$ being the depth of the network.
We use the element-wise nonlinearity $\sigma(x)=\log(\exp(x) + \exp(-x))$, which is the antiderivative of the hyperbolic tangent, i.e., $\sigma'(x)=\tanh(x)$.
For the experiments, we use $m=64$ nodes per layer and a network depth of $M=4$.
\subsection{Numerical Implementation} \label{sec:num-imp}
We tackle the control problem \cref{eq:train_loss} using the discretize-then-optimize approach. We begin by discretizing the constraints and then optimize. To do so, we first sample parameters $\bfy \sim \eta$ and then use a time discretization of $N+1$ equidistant time points $t=s_0, \ldots, s_N$ with step size $\Delta s = (T-t)/N$ to eliminate the constraints in \cref{eq:sys_dyn}. The specific discretization details for each experiment vary and are discussed in their respective sections.
This yields a state trajectory starting at $\bfz_0=\bfx$ via
\begin{equation}\label{eq:sampling}
    \bfz_{i+1} = \bfz_i + f(s_i,\bfz_{i+1};\bfy) \Delta s + g(s_i,\bfz_{i+1};\bfy) \,\bfu_i\, \Delta s, \quad i=0,\ldots,N-1,
\end{equation}
where $\bfz_i = \bfz(s_i)$ and $\bfu_i = \bfu^*(s_i,\bfz_i,\nabla_z \Phi_{\bfth} (s_i,\bfz_i;\bfy))$ is the optimal control obtained from the feedback form, that is, from~\cref{eq:opt_cont_nec_close},
computed by equating $\nabla_{\bfu}\cH=0$.
Finally, we approximate the objective functional via
\begin{equation} \label{eq:disc-obj}
    J(s_k,\bfz,\bfu;\bfy) = \Delta s \sum_{i=k}^N L(s_i,\bfz_i,\bfu_i)+ G(\bfz_N),
\end{equation}
and simultaneously, the penalty term $P_{\rm HJB, \bfy}$ in a similar manner.

In principle, any stochastic approximation approach can be used to solve the above optimization problem.
Here, we use Adam~\cite{kingma2014adam} and sample a minibatch of trajectories originating in i.i.d. samples from $\eta$ with an initial learning rate of $0.075$. During training, the learning rate gradually diminishes with an exponential decay rate of $0.975$ until it stabilizes at $0.0025$. For each training iteration, we randomly select a batch of $20$ problem parameters from the distribution $\eta$. Following the HJB framework, we evaluate controls for each problem and evolve the systems accordingly. Subsequently, we compute the losses, utilizing them to update the network weights, as illustrated in Algorithm~\ref{alg:hjb_overview}.

\begin{algorithm}[t]
\caption{Model-based training approach}\label{alg:hjb_overview}  
\begin{algorithmic}[1]
  \vspace{0.025in}
  \Require Number of training problems $\rm P$, 
  problem parameter distribution $\eta$,
  initial network weights $\bfth$, 
  Hamiltonian $H$, and loss weighting factors $(\beta_1, \beta_2, \beta_3)$
  \vspace{0.1in}  \hrule \vspace{0.1in}
  \State {\bf{Initial Setup for Training}} 
  \State {\emph{Derive Feedback Form Expression}} 
  \State $\operatorname{Feedback\_Form}: \, (-\nabla_{\bfz}\Phi_{\bfth},\, f,\, g) \, \mapsto \, \bfu$
  \State
  \State {\emph{Assemble Randomized Problems}}
  \State Sample $\rm P$ parameters \, $\{\bfy_i\}_{i=1}^{\rm P} \sim \eta$
  \State Store $\texttt{matrices}$ for $f$ and $g$ corresponding to $\{\bfy_i\}_{i=1}^{\rm P}$
  \State%
  \vspace{0.1in}  \hrule \vspace{0.1in}
  \State {\bf{HJB Training Iteration}}
  \State {\emph{Sample Batch of Problems}}
  \State  Retrieve $\texttt{matrices}$ for $\bfy_p$ with $p \sim \operatorname{Uniform}(\{1,\dots,{\rm P}\})$ 
  \State
  \State {\emph{Assess HJB Control Performance}}  
  \State  Set $s_0=0$ ,  $\bfz_0 = \bfx$ , $J=0$, and $P_{\rm{HJB}, \bfy_p} = 0$
  \While{$i \leq N$}
  \State {\emph{Derive Control from Feedback Form}}
  \State {{Evaluate network}} $\Phi_{\bfth}(s_i,\bfz_i,\bfy_p)$ and gradient $\nabla_{\bfz} \Phi_{\bfth}(s_i,\bfz_i;\bfy_p)$
  \State $\bfu_i = \operatorname{Feedback\_Form}(-\nabla_{\bfz}\Phi_{\bfth}(s_i,\bfz_i;\bfy_p), \, f,\, g)$  
  \State 
  \State {\emph{Evolve Systems in Time}}
  \State $\bfz_{i+1} \gets \operatorname{Time\_Integrator}(\bfz_i , \bfu_i ; \bfy_p)$
  \State $s_{i+1} \gets s_i + \Delta s$
  \State 
  \State {\emph{Update Intermediate Losses}}
  \State $J \gets J + L(s_i,\bfz_i,\bfu_i)  \Delta s$
  \State $P_{\rm{HJB},\bfy_p}  \gets P_{\rm{HJB},\bfy_p} + \beta_1 |H(s_i,\bfz_i, \nabla_{\bfz}\Phi_{\bfth}(s_i,\bfz_i;\bfy_p); \bfy_p)\, - \, \partial_s \Phi_{\bfth}(s_i,\bfz_i,\bfy_p)|\Delta s$%
  \EndWhile  
   \State $J \gets J + G(\bfz_N)$
   \State $P_{\rm{HJB},\bfy_p}  \gets P_{\rm{HJB},\bfy_p} +\beta_2  |G(\bfz_N - \Phi_{\bfth}(s_N, \bfz_N,\bfy_p)| +  \beta_3 |\nabla_{\bfz} G(\bfz_N) - \nabla_{\bfz} \Phi_{\bfth}(s_N, \bfz_N,\bfy_p)|$
  \State
  \State {\emph{Update Networks Weights}}
  \State $\theta \gets  \operatorname{Optimizer}(   J  \, + \,  P_{\rm{HJB},\bfy_p}  )$
\end{algorithmic}
\end{algorithm}

\section{Data-Driven Approach} \label{sec:RL}
We briefly review the data-driven approach, in particular, the Reinforcement Learning (RL) framework, with an emphasis on Actor-Critic algorithm used to approximate the optimal policy in our work;  
our presentation follows \cite{Bertsekas-RLOC,sutton2018reinforcement}.

Reinforcement learning (RL) provides a general framework for approximating policies through a series of interactions with an environment.  
More precisely, an RL environment specifies all possible system states, a description of admissible actions, along with a model governing the evolution of the system space when an action is performed. In RL language, the controls are referred to as the actions. The environment is paramount for RL algorithms, and its proper definition is crucial.

Let ${\pith}_{\bfa}$ denotes the policy, parameterized by parameters $\bfth_{\bfa}$, which is a mapping from the state space to the action space. The policy can be either deterministic, denoted as $\bfu_i = {\pith}_{\bfa}(\bfz_i)$, or stochastic, indicated as $\bfu_i \sim {\pith}_{\bfa}(\cdot \vert \bfz_t)$. Using \cref{eq:sampling}, a policy generates a sequence of states and actions, ${(\bfz_i,\bfu_i)}_{i=0}^N$, referred to as an episode.

In order to evaluate the performance of a policy, it is necessary to create functions step() and reset() governing the problem dynamics and the interaction between the agent and the environment. Aiming to address the optimization problem introduced in \cref{sec:OCF}, these functions operate as follows
\begin{itemize}[leftmargin=*]
    \item \emph{reset()}:  This function initializes or resets the environment to its initial state, $\bfz_0$, preparing it for the new episode. 
    \item \emph{step()}: Given an action $\bfu_i$ as input, determined by a policy, it returns the following:
    \begin{itemize}
        \item The next state $\bfz_{i+1}$ using the dynamics described in \cref{eq:sampling}.
        \item A scalar feedback signal representing the objective function, known as the reward,
        \begin{align*}
            r_i = \begin{cases}
                \Delta s \, L(s_i,\bfz_i, \bfu_i), &\quad i<N \\
                \Delta s \, L(s_N,\bfz_N, \bfu_N) + G(\bfz_N),  &\quad i=N
            \end{cases}.
        \end{align*}
        \item The termination status indicates the terminal time $T$ has been reached or the current episode has terminated.
    \end{itemize}
\end{itemize}
Specific environment details for the problems considered are provided in \cref{sec:num-exp}.

The goal is to identify an optimal policy in order to minimize the cumulative reward, called return\footnote{Conventional RL literature assumes that the objective function is being maximized; for the purposes of aligning with the optimal control problem setup from previous sections, we instead aim to minimize the loss (which can be interpreted as the negative reward in standard RL literature).}, given by the value function 
\begin{equation} \label{eq:return}
  \Psi(\bfth_{\bfa}) \,\,  = \,\, \E_{{\pith}_{\bfa}, \bfy\sim \eta}\left\{\, \sum_{i=0}^N r_i \right\}\,\, = \,\, \E_{{\pith}_{\bfa}, \bfy\sim \eta} \left\{\, J(s_0,\bfz,\bfu;\bfy) \right\}.
\end{equation}
When the policy is deterministic, the expectation over ${\pith}_{\bfa}$ in \cref{eq:return} is omitted, making the return equivalent to the one in the neural-HJB training defined in \cref{eq:disc-obj}.
\begin{remark}
    When starting from the initial state $\bfz_0$ and following the optimal policy ${\pith}_{\bfa}^*$, we recover the original value function $\Phi$ in \cref{eq:val_fn}.
\end{remark}

In our work, we consider the policy to be stochastic as they provide additional flexibility, robustness, and adaptability that can be advantageous in complex and uncertain environments. 

Now that we have established a loss function representative of the original optimization problem, we must also review the expression for the gradients with respect to $\bfth_{\bfa}$ for updating these parameters. This expression is provided by the following foundational result from the theory of RL.
\begin{theorem}\cite[Policy Gradient Theorem, Chapter~13]{sutton2018reinforcement}\label{thm:policy_grad}
The gradient of the cumulative reward objective function \cref{eq:return} can be expressed in terms of the policy gradients as follows: 
%
\begin{align}\label{eq:PG}
  \nabla_{\bfth_{\bfa}} \Psi(\bfth_{\bfa}) \,\, = \,\, \mathbb{E}_{{\pith}_{\bfu},\, \bfy \sim \eta}\left\{\, \sum_{i=0}^N \,   \nabla_{\bfth_{\bfa}} \log {\pith}_{\bfa}(\bfu_i | \bfz_i) {J}(s_i,\bfz,\bfu;\bfy)\right\}.
\end{align}
\end{theorem}
While the direct use of Theorem~\ref{thm:policy_grad} for policy optimization can work effectively in some cases,
the practical performance often suffers from high variance in the gradient estimate, and thus produce slow learning.
To reduce variance effects and increase stability, many policy-based RL algorithms subtract the cumulative reward by a baseline, $b(\bfs)$:
\begin{equation} \label{eq:pg_baseline}
  \nabla_{\bfth_{\bfa}} \Psi(\bfth_{\bfa}) \,\, = \,\, \mathbb{E}_{{\pith}_{\bfa},\, \bfy \sim \eta}\left\{\, \sum_{i=0}^N \,  \left(
  J(s_i,\bfz_i,\bfu_i;\bfy) - b(\bfz_i)\right)  \nabla_{\bfth_{\bfa}} \log {\pith}_{\bfa}(\bfu_i | \bfz_i) \right\}.
\end{equation}
Intuitively, making the cumulative reward smaller by subtracting it with a baseline will make smaller gradients, and thus smaller and more stable updates.

\subsection{Actor-Critic Models} \label{subsec:actor_critic}
Selecting an effective baseline poses a challenge, but the state-value function is a suitable choice, making both the value function and policy learnable. This creates a two-component RL framework known as an actor-critic model:
\begin{itemize}
\item The actor ${{\pith}_{\bfa}}$ determines actions based on the current state (i.e.\! defines the policy).
\item The critic $\Psi_{\bfth_{\bfc}}$, parameterizing the baseline with parameters ${\bfth_{\bfc}}$, assesses the approximate value of the current state (i.e.\! provides value function estimates).
\end{itemize}
and both are parameterized by neural networks.

Now that the value function is also made learnable, that means we must provide an objective for the critic network to learn the parameters $\bfth_{\bfc}$.
Since the role of the critic is to provide a baseline estimate for the return, a natural choice for this objective is the expected value of the squared error:
\begin{equation} \label{eq:critic_loss}
  P_{\rm critic}(\bfth_{\bfc}) \,\, = \,\,  \sum_{i=0}^N \, \left| \Psi_{\bfth_{\bfc}}(\bfz_i) \, - \, J(s_i,\bfz,\bfu;\bfy) \right|^2.
\end{equation}
In summary, for the actor-critic models, we are approximately solving the minimization problem given by
\begin{equation*}
    \min_{\bfth_{\bfa},\bfth_{\bfc}} \,\, \mathbb{E}_{{{\pith}_{\bfa}}, \, \bfy \sim \eta} \left\{  J(s_0,\bfz,\bfu;\bfy) + \beta_5 \,   P_{\rm critic}(\bfth_{\bfc})\right\},
\end{equation*}
where $\beta_5 \in \R_+$ controls the relative contribution of critic.
Although there are several actor-critic algorithms available, in this work we use two prominent, state-of-the-art algorithms: Proximal Policy Optimization (PPO) and Twin Delayed Deep Deterministic Policy Gradients (TD3). Both PPO and TD3 are advanced versions of the actor-critic framework discussed above, each addressing critical aspects stability, mitigating overestimation bias, or providing better exploration strategies.  For in-depth explanations and implementation specifics, we direct readers to  \cite{schulman2017proximal,fujimoto18a-TD3}.

\subsection{RL Network Architecture} \label{subsec:rl_architecture}

\begin{figure}[ht]
  \centering
  \begin{tabular}{cc}
    \includegraphics[height=1.7in]{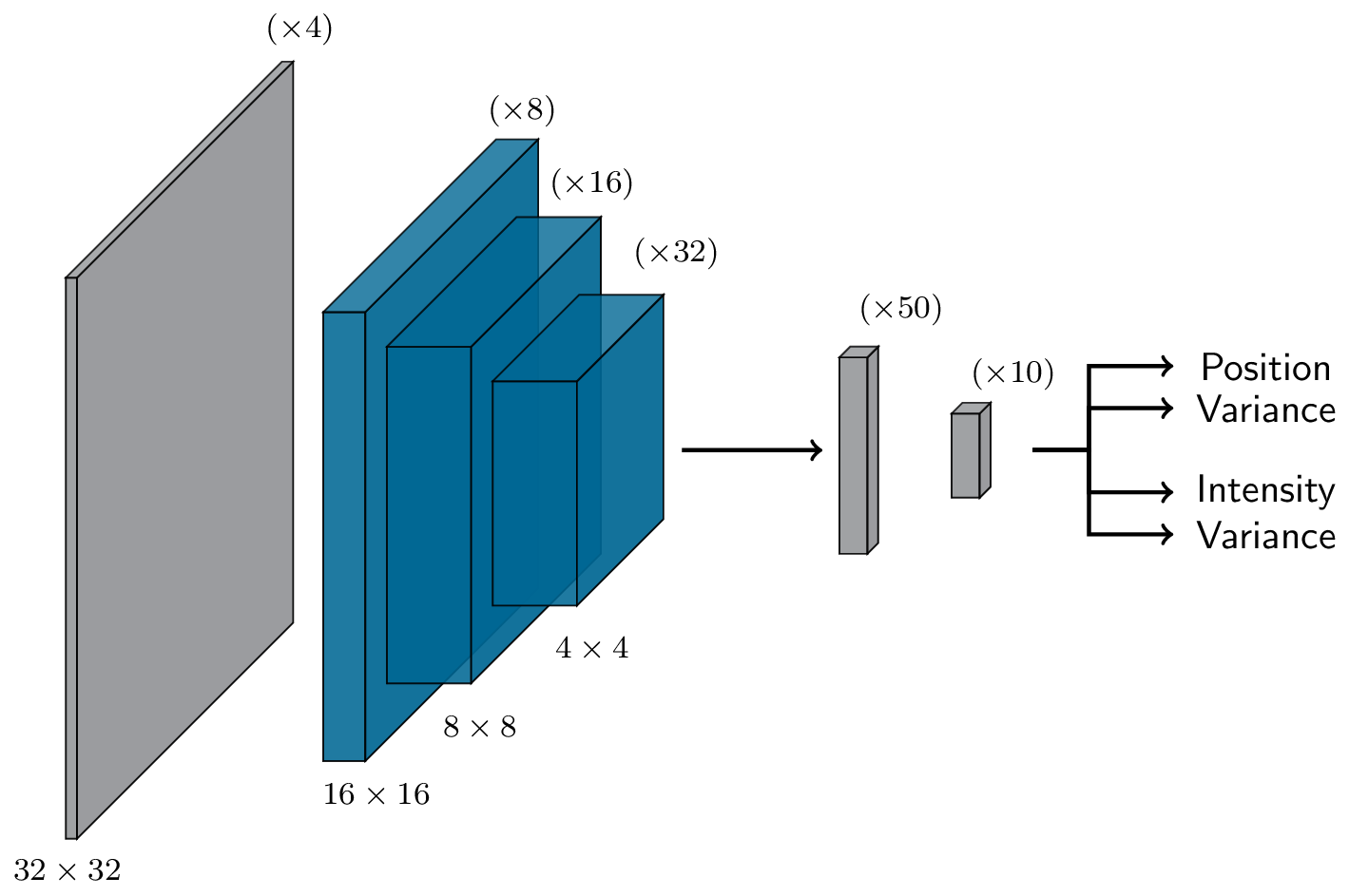} \hspace{0.1in}&
    \includegraphics[height=1.7in]{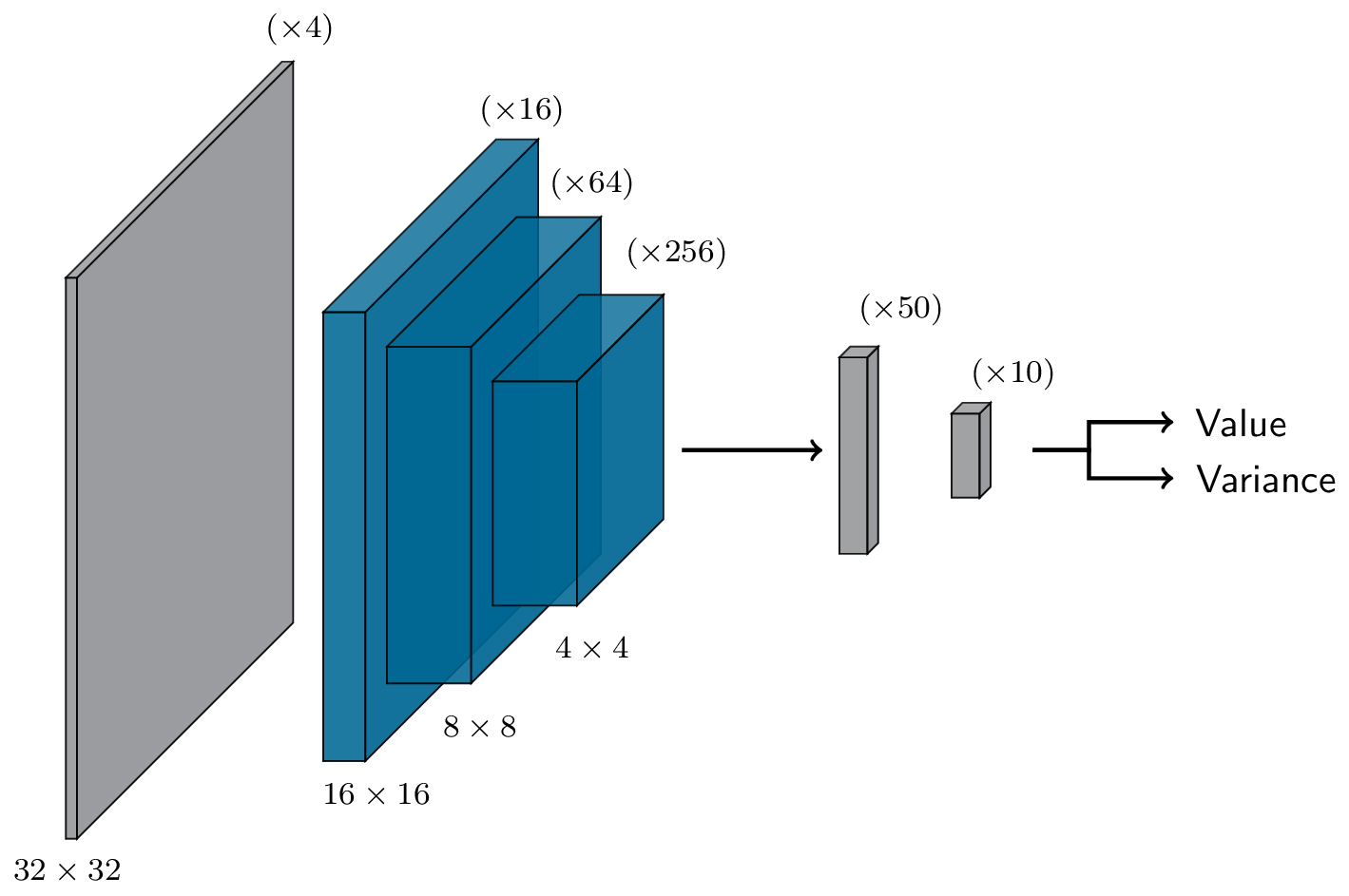}
  \end{tabular}
  \caption{Network architectures for the actor (left) and critic (right) components of the RL models.
    Both networks receive input arrays input containing the values of $\bfz(s)$, $\bfy$, and $s$. 
    Convolutional and max-pooling layers (blue) process the 
    data received from the PDE environment to extract features.  These features are then flattened and passed to dense layers (grey) 
    to form the position, variance, and value predictions.}
  \label{fig:ac_networks}
\end{figure}

Since the observation data for our experiments 
is spatially structured on a two-dimensional grid,
we employ convolutional network architectures to process the system state at each time-step.
Both the actor and critic networks share the same network architecture with the exception of the final layer which is adapted to the specific format of the actions and value outputs, respectively. 
As illustrated in \fig{fig:ac_networks}, the networks consist of three convolutional layers with $3\times 3$ kernels each followed by $2\times 2$ max-pooling layers; the series of convolutional layers is then followed by two dense layers and a final linear layer is used to produce the position/variance in the case of the actor network and the value prediction in the case of the critic network. 
All layers are equipped with hyperbolic tangent activation functions, and orthogonal initialization is used for both the convolutional and dense weights.

We train the networks using a pair of Adam optimizers, one used for the actor network and the other used for the critic, and apply an exponential decay to the learning rates which are both initially set to $10^{-4}$.
The observation states are also normalized using a exponential moving average estimate for the mean and variance of the input data in each batch.

It is essential to highlight the difference in network architectures employed for the OC and RL approaches. The OC approach necessitates an architecture with at least one continuous derivative with respect to the network inputs. %
In contrast, the RL approach only requires policies that are differentiable with respect to the network weights. Because of this, RL is compatible with most off-the-shelf architectures, %
while architecture choices for the OC approach are slightly more restricted.

\section{Numerical Results} \label{sec:num-exp}
In this section 
we compare the model-based and data-driven training approaches for a control problem involving the 
advection diffusion equation.
Since the ground-truth solutions for these problems are not available, 
we use conventional, gradient-based solutions as baselines to approximate the optimal solutions for fixed parameters in a test set.
All models were trained using 10 cores of an Intel Xeon Platinum 8176 CPU at 2.10GHz.

    In our initial experiments,
    we trained all models  using the same network architecture; however, we found the OC and RL approaches required distinct design choices to attain the best results.   When using a convolutional architecture for the OC training procedure, the gradient signals appeared to be too weak, and models failed to train from the start.  For the RL models, replacing convolutional layers with a residual network led to a noticeable drop in performance that remained after several iterations of hyperparameter tuning.  To provide a fair assessment, we present numerical results using the architecture that yields the best performance for each approach.

\subsection{Advection-Diffusion Problem Formulation} \label{subsec:problem_defs}

\begin{figure}
    \centering
    \includegraphics[width=\textwidth]{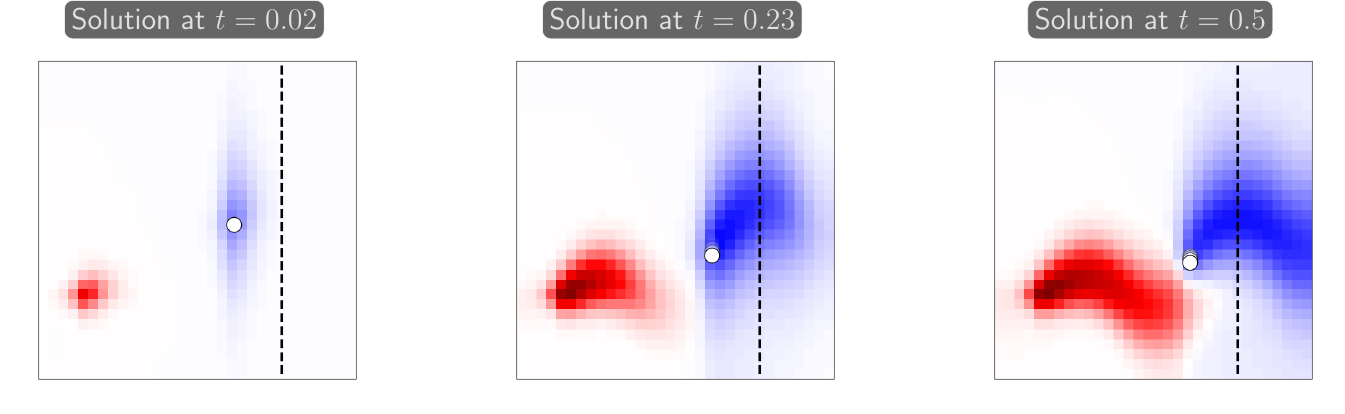}
  \caption{\label{fig:illustration_alt}Example evolution of advection diffusion system}
\end{figure}

For our numerical tests, we consider the solution $a(\cdot,\cdot)$ to an advection-diffusion PDE defined on a square domain $\Omega = [0,1]\times [0,1] \subset \mathbb{R}^2$ as part of the system state $\bfz$.
The control problem is based on a two-dimensional, continuous control vector $\bfu(s) = [u_1(s), u_2(s)]^\top$, 
aiming to prevent positive values of the solution $a(\cdot,\cdot)$ from occurring in a target region $\Omt = \{ (x_1,x_2) \in \Omega \, : \, x_1 > 0.75 \}$ at the final time $T=0.75$.
A motivating contaminant control example is illustrated in \fig{fig:illustration_alt} 
where $a$ represents contaminant concentration (in red) originating from a source location corresponding to fixed $\bfy$. The parameters $\bfy$ either contains information about the source location (Experiment 1) or both the source location and phase of advecting velocity (Experiment 2).
The goal is to prevent this contamination from reaching the target subdomain $\Omt$ which is shown to the right of the vertical dashed line. 
To achieve this, an optimal sequence of locations and magnitudes for the sink (in blue) must be found, both of 
which will depend strongly on the parameter $\bfy$ under consideration. 

For a given $\bfy \in \eta$, 
the objective function for the problem is given explicitly by:
\begin{equation}\label{eq:orig_obj_new}
  \frac{1}{2} \int_{0}^T |\bfu(s)|^2 {\dd} s + \rho \int_{\Omt} \left|\max \left\{a(T,(x_1,x_2)),0\right\} \right| \, {\dd}x_1 \, {\dd}x_2 
\end{equation}
and part of the dynamics from \cref{eq:sys_dyn} are prescribed by an advection-diffusion PDE of the form
\begin{equation}\label{eq:pde_continuous}
  {\partial_s a} + {\rm div}(\kappa\, \nabla a)  + \varepsilon_{\bfy} \, \nabla a + \varphi_{\bfy} \,\, = \,\, u_{1} \, Q(\alpha)
\end{equation}
with initial condition $a(0,x_1,x_2) = 0$ and homogeneous Neumann boundary conditions. 

The controls $\bfu$ influence system dynamics 
through the sink term `$Q$'. Here, $u_1$ determines the magnitude of the sink while $u_2$ controls the sink's velocity via the location variable $\alpha$ by:
\begin{equation*}
  \alpha(t) \,\, = \,\, \alpha(0)  \, + \, \int_0^t u_2(s) \, ds    \hspace{0.15in} \mbox{with} \hspace{0.15in} \alpha(0) = 0.5     .
\end{equation*}
For the sink term $Q$, we consider
\begin{equation*}
  Q(x_1,x_2,\alpha(t)) \,\, = \,\,   25 \cdot \exp\left(|x_1-0.6|/0.025 + |x_2-\alpha(t)|/0.15\right) .
\end{equation*}

Each parameter $\bfy \in \eta$ consists of three components $\bfy= [\bfy_{x_1}, \bfy_{x_2}, \bfy_v]$. The parameters $\bfy_{x_1}$ and $\bfy_{x_2}$ denote the $x_1$ and $x_2$ coordinates of the location of the source term, $\varphi_{\bfy}$, given by 
\begin{equation}
  \varphi_{\bfy}(x_1,x_2) \,\, = \,\, (c/\sigma_s)\, \exp\left( -(|x_1-\bfy_{x_1}| + |x_2-\bfy_{x_2}|)/\sigma_s\right)
\end{equation}
while $\bfy_v$ adjusts the phase of the velocity field, $\varepsilon_{\bfy}$.

For our first family of problems, we consider a fixed velocity field, but randomized source locations, in a horizontal setup with parameters \[\varepsilon_{\bfy}(x_1,x_2) = [25,0],\quad c=5,\quad \sigma_s=0.01,\quad \text{and}\quad  \kappa = 0.008.\] 

We then consider a relatively challenging family of problems, sinusoidal setup, which incorporates a more complex velocity field and extends problem randomization to control the phase of the velocity as well: 
\begin{equation}\label{eq:sinusoidal_velocity}
  \varepsilon_{\bfy}(x_1,x_2) \,\, = \,\, 
  \begin{bmatrix}
    (1+x_1) \cdot \sqrt{0.9^2 - (0.75\cdot\cos(1.1-x_1) \cdot \sin( 4\pi  (x_1-\bfy_{v})) )^2} \\
    \, -0.9\cdot \cos(1.1-x_1) \cdot \sin( 4\pi  (x_1-\bfy_{v}))
  \end{bmatrix}
\end{equation}
For the sinusoidal setup, we also adjust the source magnitude to $c=0.5$ and the source size to $\sigma_s=0.025$
to produce 
more natural flows for the modified velocity field.

We sample the components of $\bfy$ from the uniform distributions
\begin{equation*}
\bfy_{x_1} \sim \operatorname{Unif}(0.1,0.25)  \hspace{0.1in} , \hspace{0.15in}
\bfy_{x_2} \sim \operatorname{Unif}(0.2,0.8) \hspace{0.1in} , \hspace{0.15in}
\bfy_{v} \sim \operatorname{Unif}(-0.425,0.0)
\end{equation*}
so that the source locations are varied in the horizontal and vertical directions across the left side of the domain.
For the sinusoidal setup, the velocity field directs concentration near the source upward for values of $\bfy_v \approx -0.425$ and
directs downward for $\bfy_v \approx 0.0$.

Finally, to balance the importance of the running cost and terminal cost terms in the objective function,
we choose a scaling factor of $\rho = 40$ in \cref{eq:orig_obj_new}.  
This is generally necessary since the scale of the action space does not always match that of the state space,
and one of the loss components can dominate the other if left unweighted. 

\subsubsection{Discretization}
Since the continuous formulation results in an infinite dimensional optimal control problem,
we cannot apply numerical methods to the problem described above directly. 
Instead we consider an approximate formulation arising from the spatial discretization of \cref{eq:pde_continuous} using finite elements. 

Denoting the discretized, FEM representation of $a(s,x_1,x_2)$ by $\bfa(s)$, 
the state vector for the problem is given by: 
\begin{equation*}
  \bfz(s) \,\, = \,\, \left[ \bfa(s) \, , \, \alpha(s) \right]^\top   . 
\end{equation*}
The dynamics for the $\bfa$ component of the state vector are defined by the linear systems associated with the FEM formulation of \cref{eq:pde_continuous}. 
By construction, the state dynamics for $\alpha$ are simply $\dd_s \alpha(s)  =  u_2(s) $.

For the spatial discretization, we have selected linear,  triangular elements with nodes associated with a uniform $32\times 32$ grid.  
The FEM representation for the concentration $\bfa$ is passed to the neural networks by evaluating pointwise on
the grid locations $\{(x_1^i,x_2^j)\}_{i,j=1}^{32}$ with
$x_1^i = i\cdot {\Delta} x_1$, $x_2^j = j\cdot {\Delta} x_2$, and ${\Delta} x_1 = {\Delta} x_2 = 1/31$.
We further discretize the PDE in time via implicit Euler with a step count of $N=25$ and step size ${\Delta} t = 0.02$
to obtain the fully discretized control objective 
\begin{equation}\label{eq:orig_obj_discrete}
  J(0,\bfz,\bfu;  \bfy)  := \frac{1}{2}
  \sum_{i=0}^{25} |\bfu_i|^2 \Delta t 
  \, + \,  \rho
  \sum_{i,j=1}^{32} 1_{\Omt}(x_1^i,x_2^j)\left|\max \left\{\bfa(T),0\right\} \right|  \Delta x_1 \Delta x_2 .  
\end{equation}

\subsubsection{Feedback Form for HJB}

To apply the proposed model-based approach 
procedure, we first need to derive the feedback form for the optimal control.
This expression depends on the specific PDE under consideration and provides a natural way of incorporating system knowledge 
into the training procedure. 
As stated in \cref{eq:opt_control_nec},
the optimal control is a maximizer for the Hamiltonian evaluated at $\nabla_{\bfz}\Phi$. 
By taking the gradient of the Hamiltonian and setting the resulting expression to zero,
the feedback form for both problem setups is given by:
\begin{equation*}
\bfu(s) = \left[ -\bfQ^T \nabla_{\bfa} \Phi(s),\,  -\nabla_{\alpha} \Phi(s, \bfz(s) ; \bfy)   \right]^\top
\end{equation*}
where $\bfQ$ is the finite element matrix associated with the sink term `$Q$' in \cref{eq:pde_continuous}.

\subsubsection{Parallel Implementation of RL environments} \label{subsec:rl_parallel}

The FEM calculations used to simulate the PDE environment present a significant computational bottleneck in the training procedure.  
Fortunately, both the PPO and TD3 algorithms 
can be trained using interactions with multiple, distinct environments in parallel. 
This is possible since both algorithms decouple the actor weights used during optimization
from the network weights used in the simulations\footnote{PPO decouples weights through an `{\emph{old}}' policy, and TD3 uses a related technique using a `{\emph{target}}' policy. }.
In particular, since both algorithms update policy weights based on interactions performed by decoupled policies, 
it is possible to run multiple episodes in parallel before calculating losses and applying gradient updates. 
To take advantage of this, 
we employ OpenMPI and MPI for Python to parallelize the interactions of the actor network with several different environment realizations simultaneously.
The environments are initialized in parallel, observation data is gathered to the root node and passed to the actor/critic networks as a single batch, and the proposed actions are then broadcast back to the respective environments.
We repeat this process until all environments have completed a full episode. 
After all episodes are complete, we compute the associated returns, group transitions from all environments into minibatches, and perform gradient updates.

\subsection{Hyperparameter Selection}
Determining hyperparameters, including network architectures and optimization parameters, plays a crucial role in neural network methodologies. In our observation, the selection of hyperparameters frequently depends on the specific approach and problem at hand. Developing well-founded mathematical principles for parameter selection in these models represents a significant area for future research, but it exceeds the scope of this paper. Nonetheless, we provide our insights into hyperparameter selection for the two approaches outlined in this paper.

The proposed HJB approach requires tuning for only a small handful of parameters:
the networks' width and depth, the learning rate and decay schedule, the batch size, and three weights $(\beta_1, \beta_2, \beta_3)$
which determine how much the HJB constraints contribute to the loss.
The RL models considered in this work
have a much larger set of hyperparameters and require a substantial number of system solves to
identify the correct training settings. 
This is particularly true of scaling and clipping parameters that are used to prevent the RL training procedure from becoming unstable. 
When these parameters are not set correctly, %
the models' performance will often deteriorate midway through training, and the procedure must be restarted.

\subsection{Baseline via Gradient-Based Solutions}
\label{sec:baseline}

To understand how the proposed, amortized procedure compares with more conventional approaches,
we now review how gradient-based methods handle this class of control problems.
These methods apply to a single problem realization (i.e., fixed value of $\bfy$)
and rely on gradient information, namely $\nabla_{\bfu} J$, to identify an optimal control.
For this work, we have implemented a gradient-based solver
using FEM matrices extracted from FEniCS
and automatic differentiation provided by PyTorch to retrieve the necessary gradients. 
We then use the 
limited memory Broyden-Fletcher-Goldfarb-Shanno algorithm
to search the control space for an optimal solution.  
We use the solutions obtained from this gradient-based search as proxies for the true optimal solutions (which are intractable for the problems considered in this work).

Here, it is important to note that the control obtained from the gradient-based method is specific to the given parameter. 
While the control may perform well for nearby parameters, in general, changes in the parameter require another solve. 
In many cases, the original control 
does not even provide an effective initialization point for new problems (since the optimal control can change considerably when system dynamics are modified). 
Because of this, the initialization procedure is another important and non-trivial aspect of gradient-based methods.
Even though these methods are designed to solve a single problem, unlike amortized approaches that target a number of problems at once, 
the gradient-based search can fail to converge to an optimal control without proper initialization.

\subsection{Assessment of HJB and RL Performance} \label{subsec:results}

\begin{figure}
    \centering
    \includegraphics[width=\textwidth]{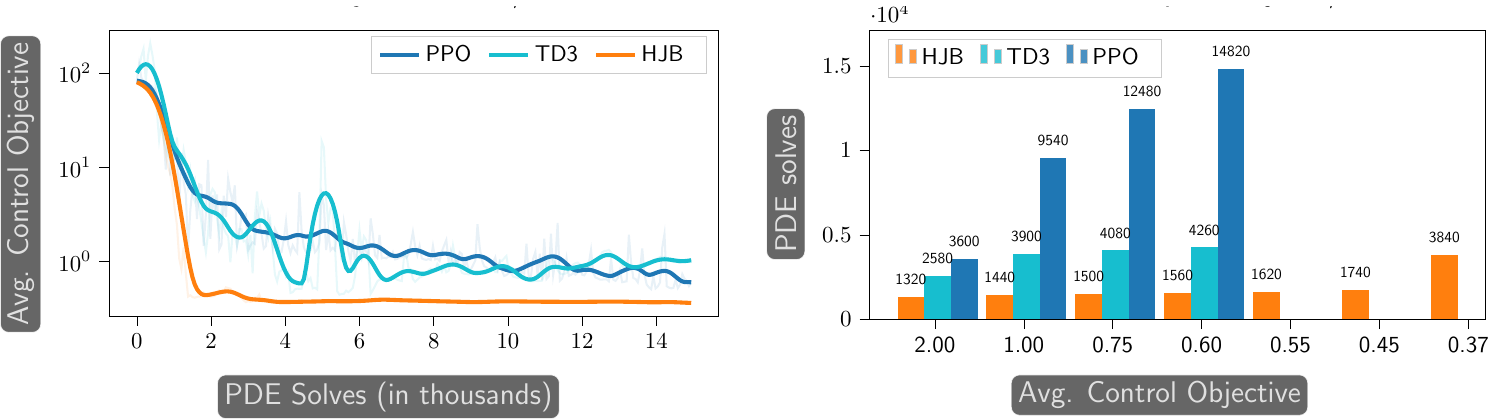}
    \caption{Horizontal problem setup: (left) Validation loss during training and (right) number of PDE solves required for different target accuracies of control objective.} 
    \label{fig:horz_dual}
\end{figure}

For the first set of experiments, 
we consider the horizontal velocity setup with two controls 
corresponding to the magnitude and velocity of the sink term.
The RL and HJB models are tasked with positioning the sink in locations
that block the 
concentration from entering 
the target region while using
as little movement as possible. 
To evaluate each model's performance,
we assess the accuracy of the models on a set of fixed validation problems throughout training.
For the validation set,
we have selected parameter values 
$\bfy_{x_1} \in \{0.125 , 0.225\}$ and $\bfy_{x_2} \in \{0.25 , 0.4, 0.5, 0.6 , 0.75\}$.
This corresponds to $10$ realizations of the horizontal setup
with source locations spread across 
the left side of the domain. 
The results of this experiment are summarized in \fig{fig:horz_dual} 
where we observe that the HJB framework 
yields several, noticeable improvements
over both 
RL models.

The most clear, and practically significant, distinction between the HJB and RL training
procedures is the rapid drop in the HJB validation loss at the very start of training.
This is perhaps to be expected thanks to the gradient information exposed by the HJB
training formulation, while the RL models are initially reliant on an uninformed, trial-and-error
state of exploration.
For practical applications, this is a key feature  
of the HJB training procedure;
the model is able to identify effective control policies using far fewer forward-model PDE solves
than its RL counterparts.
By relying on fewer forward solves, the HJB approach has the potential to be applied to more
computationally expensive simulations where performing the tens of thousands of queries required by RL models is simply impractical. 

To asses how well the model-based approach 
framework extends
to more complicated systems,
we now turn to evaluating its performance on the sinusoidal setup. 
This setup introduces more complex system dynamics 
which incorporate 
non-uniform velocity fields for advection. 
We also add another layer 
of randomization 
by varying the phase of the velocity field
according to 
an additional problem parameter $\bfy_{v}$. 
These variations in the velocity field are 
designed to reflect uncertain ambient conditions (such as weather) which
may need to be taken into account when developing control strategies for practical applications. 
This form of randomization influences how systems evolve in time 
and will help us gauge how well the model-based approach and RL approaches
are able to 
amortize over problems with
substantially different dynamics.

For the sinusoidal setup, we define the validation set 
using parameter values selected from
$\bfy_{x_1} \in \{0.125 , 0.225\}$, $\bfy_{x_2} \in \{0.25 , 0.4, 0.5, 0.6 , 0.75\}$, and $\bfy_v \in \{-0.35,-0.2125,-0.1\}$.
This yields $30$ problem realizations with source locations spread across the left side of the domain
and velocity fields which direct concentration near these source terms upward, horizontally, as well as downward.

\begin{figure}
    \centering
    \includegraphics[width=\textwidth]{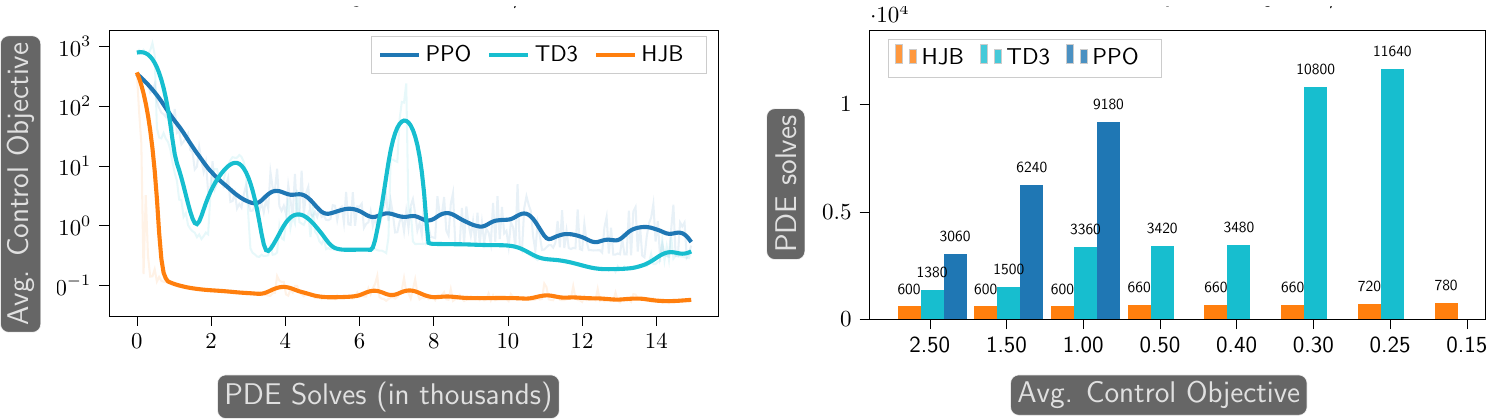}
    \caption{Sinusoidal problem setup: (left) Validation loss during training and (right) number of PDE solves required for different target accuracies of control objective.} 
    \label{fig:sine_dual}
\end{figure}

The results for the sinusoidal setup are depicted in \Cref{fig:sine_dual}, where we once again observe a clear advantage in data efficiency for the HJB model. The HJB model achieves an average control objective of $0.15$ using only 780 PDE solves. However, both RL models fail to achieve an average control objective below $0.25$ even with the full 15,000 solves.

\begin{figure}
    \centering
    \includegraphics[width=\textwidth]{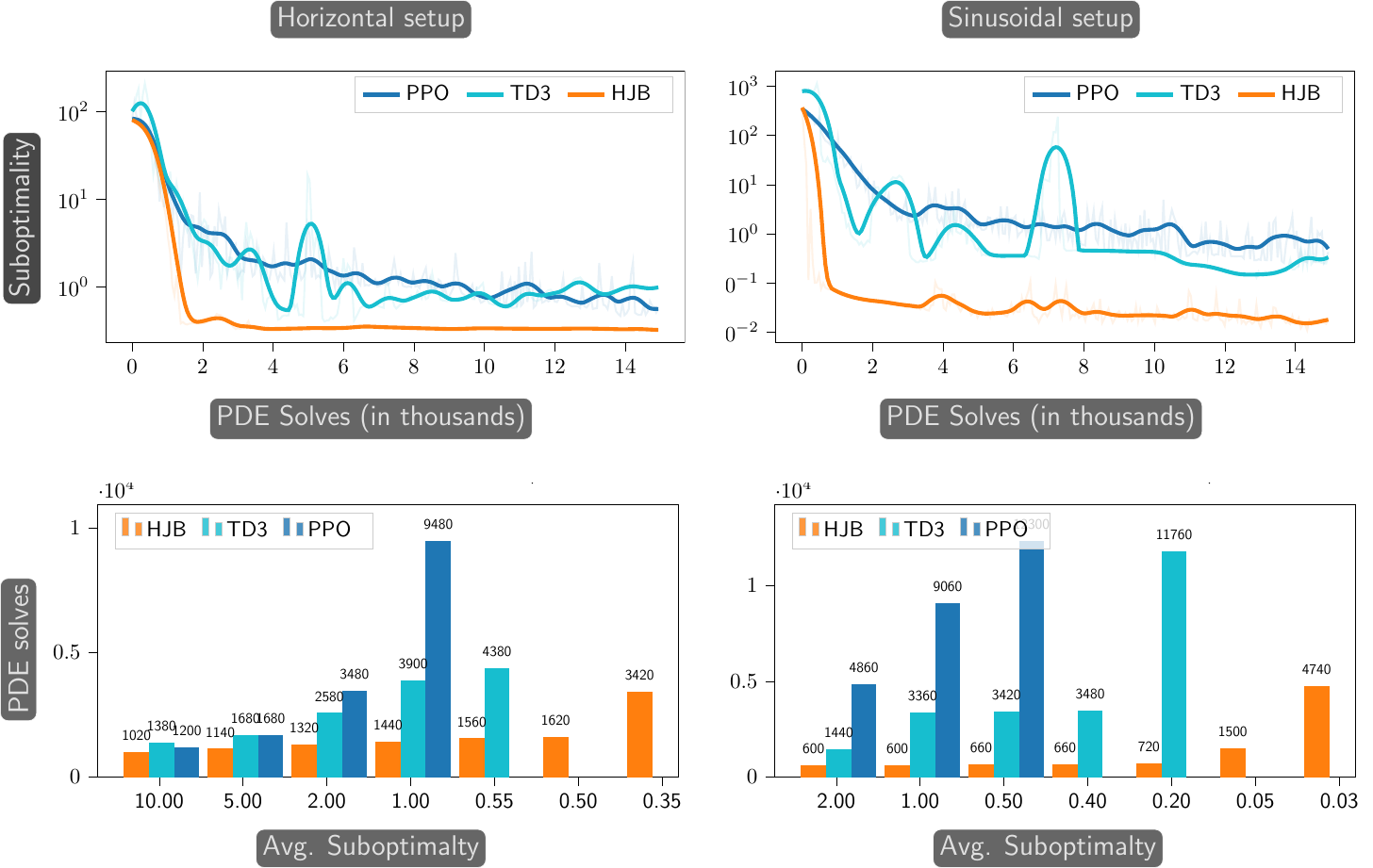}
    \caption{Suboptimality, relative to the baseline, on validation problems for the horizontal (left column) and sinusoidal (right column) problem setups.}
    \label{fig:subopt_plots}
\end{figure}

In \Cref{fig:subopt_plots}, we also present the suboptimality of HJB and RL approaches compared to the baseline, as discussed in \Cref{sec:baseline}. It can be observed that the HJB model achieves much lower suboptimality than RL models with significantly fewer PDE solves. This difference is particularly pronounced in the sinusoidal setup, where the HJB approach can achieve a suboptimality of $0.03$ in 4740 PDE solves, while RL approaches fail to go below $0.20$.

This faster convergence can be partially attributed to the training algorithms used in RL. While both RL models succeed in preventing positive concentrations from reaching the target region, they do so with unnecessary movement. The RL policies often exhibit bias in initial movement in one direction (e.g., drifting downward at the start of a simulation) and then correct with a sharp reversal in direction if needed. These inefficiencies result in much larger running costs for the RL models and prevent them from reaching the control objective levels attained by the HJB approach.

\section{Discussion}
We consider deterministic, finite-time optimal control problems with unknown or uncertain parameters and aim to approximate optimal policies for all parameters with neural networks.
We compare to training paradises.
On the one hand, our model-based approach parameterizes the value function using a neural network and derives the policy from the feedback form. By leveraging control theory, this approach incorporates the system's physics by penalizing the Hamilton-Jacobi-Bellman equation satisfied by the value function during training. 
On the other hand, the data-driven approach adopts state-of-the-art actor-critic RL frameworks and trains separate networks for the policy (actor) and the value function (critic). Both networks are trained exclusively from observational data. Despite drawing inspiration from control theory, particularly Dynamic Programming, both approaches exhibit distinct behaviors in practical applications.

Although we strived to use almost identical learning problems in both paradigms, we had to allow for some differences to account for their assumptions.
Most importantly, our model-based approach had to use a simpler network architecture to ensure the differentiability with respect to the inputs needed to apply the feedback form. While non-differentiable off-the-shelf CNN architectures worked well for the data-driven approach, the training did not yield adequate results using the simpler architecture.

We apply both approaches to a model problem inspired by contaminant containment to assess accuracy and efficiency. In this scenario, the state of the system is governed by a two-dimensional advection-diffusion partial differential equation (PDE), and the objective is to prevent the contaminant from reaching a target region by controlling the location and intensity of a sink, with the parameter reflecting various scenarios for the source locations and velocity fields.

While both approaches follow an unsupervised training regime, our numerical results reveal that the model-based approach approximates a reasonable policy with much lower sub-optimality, when compared to a gradient-based baseline. Notably, it achieves this with far fewer PDE solves than RL~(see \Cref{fig:subopt_plots}). The superior performance of the model-based approach can be attributed to the incorporation of physics and system derivatives guided by control theory. In contrast, lacking awareness of this information, RL demands a larger number of PDE solutions. During RL training, many models showed promising learning for several thousand PDE solves but eventually diverged and were unable to recover. This phenomenon never occurred in the model-based approach.

Despite RL's data-prohibitive nature, its non-intrusive and gradient-free characteristics make it widely applicable and suitable for experimenting with new problems lacking well-defined models or facing challenges in efficiently computing feedback forms and system derivatives. As part of future work, we aim to extend our model-based approach to such scenarios and revisit the comparison against RL.

Another advantage of the model-based approach lies in its simpler hyperparameter setup, requiring fewer than half the number of hyperparameters compared to RL approaches. Effective hyperparameters identified for one problem setup could be seamlessly applied to other problems without modification. Conversely, RL training procedures proved sensitive to hyperparameter choices, often necessitating re-tuning between problems for convergence. This again emphasizes the data efficiency of the model-based approach.

As a part of future work, we plan to compare the methods for more complex problems and investigate their scalability using increasingly finer discretizations.

\section*{Acknowledgments}

LR's and DV's work was supported in part by NSF awards DMS 1751636 and DMS 2038118, AFOSR grant FA9550-20-1-0372, and US DOE Office of Advanced Scientific Computing Research
  Field Work Proposal 20-023231.
NW's and BV's work was supported  in part by Sandia National Laboratories which is a multimission laboratory managed and operated by National Technology and Engineering Solutions of Sandia LLC, a wholly owned subsidiary of Honeywell International, Inc., for the U.S. Department of Energy's National Nuclear Security Administration under Contract DE-NA0003525, and US Department of Energy, and Office of Advanced Scientific Computing Research, Field Work Proposal 20-023231.
\bibliographystyle{plain}
\bibliography{ref}%

\begin{thebibliography}{10}

\bibitem{amin2021survey}
Susan Amin, Maziar Gomrokchi, Harsh Satija, Herke van Hoof, and Doina Precup.
\newblock A survey of exploration methods in reinforcement learning.
\newblock {\em arXiv preprint arXiv:2109.00157}, 2021.

\bibitem{Ballarin2017NumericalMO}
Francesco Ballarin, Elena Faggiano, Andrea Manzoni, Alfio~Maria Quarteroni,
  Gianluigi Rozza, Sonia Ippolito, Carlo Antona, and Roberto Scrofani.
\newblock Numerical modeling of hemodynamics scenarios of patient-specific
  coronary artery bypass grafts.
\newblock {\em Biomechanics and Modeling in Mechanobiology}, 16:1373--1399,
  2017.

\bibitem{bansal2021deepreach}
Somil Bansal and Claire~J. Tomlin.
\newblock {DeepReach}: A deep learning approach to high-dimensional
  reachability.
\newblock In {\em IEEE International Conference on Robotics and Automation
  (ICRA)}, pages 1817--1824, 2021.

\bibitem{Berner2019Dota2W}
Christopher Berner, Greg Brockman, Brooke Chan, Vicki Cheung, Przemyslaw
  Debiak, Christy Dennison, David Farhi, Quirin Fischer, Shariq Hashme,
  Christopher Hesse, Rafal J{\'o}zefowicz, Scott Gray, Catherine Olsson,
  Jakub~W. Pachocki, Michael Petrov, Henrique~Pond{\'e} de~Oliveira~Pinto,
  Jonathan Raiman, Tim Salimans, Jeremy Schlatter, Jonas Schneider, Szymon
  Sidor, Ilya Sutskever, Jie Tang, Filip Wolski, and Susan Zhang.
\newblock Dota 2 with large scale deep reinforcement learning.
\newblock {\em ArXiv}, abs/1912.06680, 2019.

\bibitem{Bertsekas-RLOC}
Dimitri~P. Bertsekas.
\newblock {\em Reinforcement learning and optimal control}.
\newblock Athena Scientific Optimization and Computation Series. Athena
  Scientific, Belmont, MA, [2019] \copyright 2019.
\newblock Second printing with editorial revisions.

\bibitem{traffic-flow-2010}
Rodrigo~C. Carlson, Ioannis Papamichail, Markos Papageorgiou, and Albert
  Messmer.
\newblock Optimal motorway traffic flow control involving variable speed limits
  and ramp metering.
\newblock {\em Transportation Science}, 44(2):238--253, 2010.

\bibitem{Chen2018TopKOC}
Minmin Chen, Alex Beutel, Paul Covington, Sagar Jain, Francois Belletti, and
  Ed~H. Chi.
\newblock Top-k off-policy correction for a reinforce recommender system.
\newblock {\em Proceedings of the Twelfth ACM International Conference on Web
  Search and Data Mining}, 2018.

\bibitem{demo2021extended}
Nicola Demo, Maria Strazzullo, and Gianluigi Rozza.
\newblock {A}n extended physics informed neural network for preliminary
  analysis of parametric optimal control problems.
\newblock {\em arXiv:2110.13530}, 2021.

\bibitem{flemingsoner06}
Wendell~H. Fleming and H.~Mete Soner.
\newblock {\em Controlled {M}arkov Processes and Viscosity Solutions},
  volume~25 of {\em Stochastic Modelling and Applied Probability}.
\newblock Springer, New York, second edition, 2006.

\bibitem{fujimoto18a-TD3}
Scott Fujimoto, Herke van Hoof, and David Meger.
\newblock Addressing function approximation error in actor-critic methods.
\newblock In Jennifer Dy and Andreas Krause, editors, {\em Proceedings of the
  35th International Conference on Machine Learning}, volume~80 of {\em
  Proceedings of Machine Learning Research}, pages 1587--1596. PMLR, 10--15 Jul
  2018.

\bibitem{He_2016_CVPR}
Kaiming He, Xiangyu Zhang, Shaoqing Ren, and Jian Sun.
\newblock Deep residual learning for image recognition.
\newblock In {\em Proceedings of the IEEE Conference on Computer Vision and
  Pattern Recognition (CVPR)}, June 2016.

\bibitem{henderson2018deep}
Peter Henderson, Riashat Islam, Philip Bachman, Joelle Pineau, Doina Precup,
  and David Meger.
\newblock Deep reinforcement learning that matters.
\newblock In {\em Proceedings of the AAAI conference on artificial
  intelligence}, volume~32, 2018.

\bibitem{hwang2022solving}
Rakhoon Hwang, Jae~Yong Lee, Jin~Young Shin, and Hyung~Ju Hwang.
\newblock Solving {PDE}-constrained control problems using operator learning.
\newblock In {\em Proceedings of the AAAI Conference on Artificial
  Intelligence}, volume~36, pages 4504--4512, 2022.

\bibitem{Hwangbo2019LearningAA}
Jemin Hwangbo, Joonho Lee, Alexey Dosovitskiy, Dario Bellicoso, Vassilios
  Tsounis, Vladlen Koltun, and Marco Hutter.
\newblock Learning agile and dynamic motor skills for legged robots.
\newblock {\em Science Robotics}, 4, 2019.

\bibitem{rlblogpost}
Alex Irpan.
\newblock Deep reinforcement learning doesn't work yet.
\newblock 2018.

\bibitem{kang2019algorithms}
Wei Kang, Qi~Gong, and Tenavi Nakamura-Zimmerer.
\newblock Algorithms of data development for deep learning and feedback design.
\newblock {\em arXiv:1912.00492}, 2019.

\bibitem{kingma2014adam}
Diederik~P Kingma and Jimmy Ba.
\newblock Adam: A method for stochastic optimization.
\newblock {\em arXiv:1412.6980}, 2014.

\bibitem{kunisch2020semiglobal}
Karl Kunisch and Daniel Walter.
\newblock Semiglobal optimal feedback stabilization of autonomous systems via
  deep neural network approximation.
\newblock {\em ESAIM: Control, Optimisation and Calculus of Variations}, 27,
  2021.

\bibitem{kunisch2023optimal}
Karl Kunisch and Daniel Walter.
\newblock Optimal feedback control of dynamical systems via value-function
  approximation.
\newblock {\em arXiv:2302.13122}, 2023.

\bibitem{Li-Verma-Ruthotto-2022}
Xingjian Li, Deepanshu Verma, and Lars Ruthotto.
\newblock A neural network approach for stochastic optimal control.
\newblock {\em arXiv:2209.13104}, 2022.

\bibitem{lopez2019solutions}
Victor~G Lopez, Frank~L Lewis, Yan Wan, Edgar~N Sanchez, and Lingling Fan.
\newblock Solutions for multiagent pursuit-evasion games on communication
  graphs: Finite-time capture and asymptotic behaviors.
\newblock {\em IEEE Transactions on Automatic Control (TAC)}, 65(5):1911--1923,
  2019.

\bibitem{lu2021physics}
Lu~Lu, Raphael Pestourie, Wenjie Yao, Zhicheng Wang, Francesc Verdugo, and
  Steven~G Johnson.
\newblock Physics-informed neural networks with hard constraints for inverse
  design.
\newblock {\em SIAM Journal on Scientific Computing}, 43(6):B1105--B1132, 2021.

\bibitem{lye2021iterative}
Kjetil~O Lye, Siddhartha Mishra, Deep Ray, and Praveen Chandrashekar.
\newblock {I}terative surrogate model optimization ({ISMO}): an active learning
  algorithm for {PDE} constrained optimization with deep neural networks.
\newblock {\em Computer Methods in Applied Mechanics and Engineering},
  374:113575, 2021.

\bibitem{Mnih2013PlayingAW}
Volodymyr Mnih, Koray Kavukcuoglu, David Silver, Alex Graves, Ioannis
  Antonoglou, Daan Wierstra, and Martin~A. Riedmiller.
\newblock Playing atari with deep reinforcement learning.
\newblock {\em arXiv}, 1312.5602, 2013.

\bibitem{mowlavi2022optimal}
Saviz Mowlavi and Saleh Nabi.
\newblock {O}ptimal control of {PDE}s using physics-informed neural networks.
\newblock {\em Journal of Computational Physics}, page 111731, 2022.

\bibitem{nikishin2018improving}
Evgenii Nikishin, Pavel Izmailov, Ben Athiwaratkun, Dmitrii Podoprikhin, Timur
  Garipov, Pavel Shvechikov, Dmitry Vetrov, and Andrew~Gordon Wilson.
\newblock Improving stability in deep reinforcement learning with weight
  averaging.
\newblock In {\em Uncertainty in artificial intelligence workshop on
  uncertainty in Deep learning}, 2018.

\bibitem{onken2020neural}
Derek Onken, Levon Nurbekyan, Xingjian Li, Samy~Wu Fung, Stanley Osher, and
  Lars Ruthotto.
\newblock A neural network approach for high-dimensional optimal control.
\newblock {\em arXiv:2104.03270}, 2021.

\bibitem{rozza2012reduction}
Gianluigi Rozza, Andrea Manzoni, and Federico Negri.
\newblock {R}eduction strategies for {PDE}-constrained optimization problems in
  haemodynamics.
\newblock In {\em Proceedings of the 6th European Congress on Computational
  Methods in Applied Sciences and Engineering}, number CONF, pages 1748--1769.
  Vienna Technical University, 2012.

\bibitem{schulman2017proximal}
John Schulman, Filip Wolski, Prafulla Dhariwal, Alec Radford, and Oleg Klimov.
\newblock Proximal policy optimization algorithms.
\newblock {\em arXiv:1707.06347}, 2017.

\bibitem{Silver2016MasteringTG}
David Silver, Aja Huang, Chris~J. Maddison, Arthur Guez, L.~Sifre, George
  van~den Driessche, Julian Schrittwieser, Ioannis Antonoglou, Vedavyas
  Panneershelvam, Marc Lanctot, Sander Dieleman, Dominik Grewe, John Nham, Nal
  Kalchbrenner, Ilya Sutskever, Timothy~P. Lillicrap, Madeleine Leach, Koray
  Kavukcuoglu, Thore Graepel, and Demis Hassabis.
\newblock Mastering the game of go with deep neural networks and tree search.
\newblock {\em Nature}, 529:484--489, 2016.

\bibitem{Silver2018AGR}
David Silver, Thomas Hubert, Julian Schrittwieser, Ioannis Antonoglou, Matthew
  Lai, Arthur Guez, Marc Lanctot, L.~Sifre, Dharshan Kumaran, Thore Graepel,
  Timothy~P. Lillicrap, Karen Simonyan, and Demis Hassabis.
\newblock A general reinforcement learning algorithm that masters chess, shogi,
  and go through self-play.
\newblock {\em Science}, 362:1140 -- 1144, 2018.

\bibitem{strazzullo2018model}
Maria Strazzullo, Francesco Ballarin, Renzo Mosetti, and Gianluigi Rozza.
\newblock Model reduction for parametrized optimal control problems in
  environmental marine sciences and engineering.
\newblock {\em SIAM Journal on Scientific Computing}, 40(4):B1055--B1079, 2018.

\bibitem{sutton2018reinforcement}
Richard~S Sutton and Andrew~G Barto.
\newblock {\em Reinforcement learning: An introduction}.
\newblock MIT press, 2018.

\bibitem{Trltzsch2010OptimalCO}
Fredi Tr{\"o}ltzsch.
\newblock {\em Optimal Control of Partial Differential Equations: Theory,
  Methods and Applications}, volume 112 of {\em Graduate Studies in
  Mathematics}.
\newblock 2010.

\bibitem{wang2021fast}
Sifan Wang, Mohamed~Aziz Bhouri, and Paris Perdikaris.
\newblock Fast {PDE}-constrained optimization via self-supervised operator
  learning.
\newblock {\em arXiv:2110.13297}, 2021.

\bibitem{Wu2017LearningTS}
Jiajun Wu, Erika Lu, Pushmeet Kohli, Bill Freeman, and Joshua~B. Tenenbaum.
\newblock Learning to see physics via visual de-animation.
\newblock In {\em NIPS}, 2017.

\bibitem{xu2021machine}
Mengfei Xu, Shufang Song, Xuxiang Sun, Wengang Chen, and Weiwei Zhang.
\newblock {M}achine learning for adjoint vector in aerodynamic shape
  optimization.
\newblock {\em Acta Mechanica Sinica}, pages 1--17, 2021.

\bibitem{Yin2023AONNAA}
Peng-Heng Yin, Guangqiang Xiao, Keju Tang, and Chao Yang.
\newblock {AONN}: An adjoint-oriented neural network method for all-at-once
  solutions of parametric optimal control problems.
\newblock {\em arXiv:2302.02076}, 2023.

\bibitem{yong_zhou-1999}
Jiongmin Yong and Xun~Yu Zhou.
\newblock {\em Stochastic controls}, volume~43 of {\em Applications of
  Mathematics (New York)}.
\newblock Springer-Verlag, New York, 1999.
\newblock Hamiltonian systems and HJB equations.

\bibitem{Zhou_2021}
Mo~Zhou, Jiequn Han, and Jianfeng Lu.
\newblock Actor-critic method for high dimensional static
  {H}amilton--{J}acobi--{B}ellman partial differential equations based on
  neural networks.
\newblock {\em {SIAM} Journal on Scientific Computing}, 43(6):A4043--A4066, jan
  2021.

\bibitem{zhou2023policy}
Mo~Zhou and Jianfeng Lu.
\newblock A policy gradient framework for stochastic optimal control problems
  with global convergence guarantee.
\newblock {\em arXiv:2302.05816}, 2023.

\bibitem{deepRL-chapter7-2020}
Hao~Dong Zihan~Ding.
\newblock Challenges of reinforcement learning.
\newblock In Shanghang~Zhang Hao~Dong, Zihan~Ding, editor, {\em Deep
  Reinforcement Learning: Fundamentals, Research, and Applications}, chapter~7,
  pages 249--272. Springer Nature, 2020.
\newblock \url{http://www.deepreinforcementlearningbook.org}.

\end{thebibliography}
\end{document}